# LEARNING OF THE STUDENTS IN A REPRODUCTION OF FIGURE BY FOLDING

Claire Guille-Biel Winder ESPE de Nice

*We are interested in the learning of 6 to 7 years old children in implementations of a situation of reproduction of figure by folding presented in the first part of this article. In the second part we expose our problem as well as our hypothesis and our methodology. The third part presents the mathematical and didactical potentialities of the PLIOX situation within the framework of the Theory of Didactic Situations of Brousseau, in support of Berthelot and Salin's works, and according to a cognitive and semiotic point of view of Duval. According to this previous part, we clarify and analyze the results obtained from implementations in two classrooms (part four).*

## PRESENTATION OF THE PLIOX SITUATION

We designed by "PLIOX"[1] a squared paper separated in four squared colored zones (yellow, green, blue and red) as following:

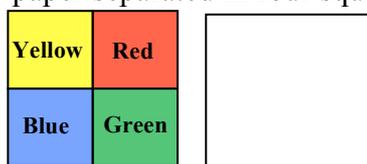

Figure 1: Recto and verso of the PLIOX

The generic PLIOX situation consists in the following didactic situation (according to Brousseau sense): reproduction of a model figure by folding a PLIOX. The authorized folds are presented on Figure 2 and correspond to the symmetric axes of square (diagonals and medians) and to the symmetric axes (diagonals and medians) of the so-called "secondary squares"(*i.e* the four colored squares).

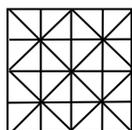

Figure 2 : Authorized folds

The PLIOX situation can be seen as an "origami" in the sensitive space, that is to say it corresponds to a problem in concrete mode (Duval, Godin & Perrin-Glorian, 2005), and also to a spatial problem. The corresponding space is the 3D micro-space (Berthelot & Salin, 1992): the PLIOX is accessible to vision; turning, returning, folding and unfolding the PLIOX are the possible manipulations (3D gestures with tactile and visual coordination); all the movements of the PLIOX in regard to the subject are possible, as well as those of the subject; the subject is outside the considered space, which contains the object.

Moreover, the PLIOX situation can be considered according to a geometrical point of view. The three dimensional object PLIOX is indeed very flat so that its form can be related to a representation of a two dimensional geometrical figure (possibly decomposed into other sub-figures). In primary school, we believe that a figure is a drawing whose properties can be specified within the theoretical geometrical framework by considering it as a representation of a geometrical

---

[1] The geometrical shape is then dissociated of the artifact (Rabardel,1995): a PLIOX is an artifact that becomes an instrument as soon as it is used to reproduce a figure by folding.

figure (Duval, Godin & Perrin-Glorian, 2005). In this case, a figure can be obtained by juxtaposition or superposition of shapes (Perrin-Glorian, Mathé & Leclerc, 2013). Consequently, we consider the reproduction of figures by folding a PLIOX as geometrical problems.

**PROBLEM, HYPOTHESIS AND METHODOLOGY**

In order to look into the learning of students in the PLIOX situation, we are interested in the mathematical and didactical potentialities of the situation and in the mathematical knowledge that circulate during implementation. Moreover, according to Duval, Godin and Perrin-Glorian (2005), we make the hypothesis that language can reveal some of these knowledge as well as their evolution. We study the potentialities of the situation through an a priori analysis and the language through the transcription of the recording of two implementations realized into classrooms of first grade primary school (6 to7 years old children). In each classroom, the progress of the situation consists of three sessions[2]: the first one is devoted to the manufacturing of the artifact by students and the two others to the reproduction of several figures by folding, according to the scheme:

- Phase 1: Presentation of the model figure in the board.
- Phase 2: Short collective analysis of the model figure.
- Phase 3: Individual reproduction of the model figure by folding.
- Phase 4: Collective highlighting of procedures and validation.
- Phase 5: Mathematical synthesis.

**ELEMENTS OF *A PRIORI* ANALYSIS**

In this part, we resume the results presented in details in (Guille-Biel Winder, 2014). The PLIOX situation involves spatial knowledge, because of the orientation and the positions of figures that appear into model figures, and also geometrical knowledge (mainly about the square, the rectangle, the right-angle triangle and their relationships, and the specific straight lines (diagonals and medians)). Moreover, two points of view on a same figure can be involved in this situation: first, what is immediately recognized (the colored zones and the global form according to Gestalt Theory) and secondly the figures who can be identified after taking account of certain folds. For this case, knowledge and recognition of elementary two-dimensional figural units require perceptual apprehension and operative (mereologic, optic and place) apprehension (Duval, 1999).

For example, reproducing Figure 3 requires to identify the external figure (the isosceles right-angle triangle) and to take care of the three colored zones (green, blue and yellow) and their relative positions (but not necessarily to recognize them as triangles or square). A change of orientation is also necessary if the model is presented as in Figure 3. This model highlights a relationship between square and isosceles right-angle triangle: dividing a square into two parts can lead to two isosceles right-angle triangles. The realization of the reproduction requires folding according to a diagonal of the square. Then the authorized folds allow recognizing two blue triangles, a big green and blue triangle, but also a blue and green square or a green blue and yellow rectangle, …

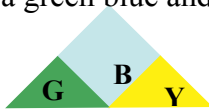

Figure 3: A model figure

---

[2] Analysis do not take account of the first sessions in the classrooms.

Still the role and the status of knowledge within the various phases can be identified. In phase 3 (action) the spatial knowledge remain implicit; some skills about axial symmetry (mostly into the square) are used; the control of the action (by visual perception) take implicitly into account the nature of various figures that exist into the model figure (geometrical knowledge), their orientation, the position and the orientation of the colored zones (spatial knowledge). Formulation phases (phases 2 and 4) allow to explicit the mereologic decomposition induced by color as well as the spatial organization of these figural units; the request to justify the identification of a figure leads gradually the children to become aware of some properties (the number of sides and summits of the square, of the rectangle, …); some relationships between some of the shapes can also been highlighted. The attempted formulations in phases 2 and 4 describe strategies and gestures in the three-dimensional space (that can be helpful for some children to achieve the task) and explicit the spatial knowledge used in phase 3. During these phases, the spatial vocabulary (over, down, behind, …) and the geometric vocabulary (square, triangle, diagonal, point, …) appear as necessary to the communication. Phase 5 results on a local institutionalization about vocabulary, geometric objects and their properties, or highlighted relationships.

**ANALYSIS OF THE LANGUAGE IN THE IMPLEMENTATIONS**

We first identify the words used by students and teacher to mean geometric objects: the nomination of usual shapes uses geometric vocabulary (square, triangle, …) or refers to the common language (color or else); the nomination of rights and segments uses geometric vocabulary (right, side, …), or refers to the common language ("line", "trace", …) or to the sensitive space ("fold", "edge", …); the nomination of points uses geometric vocabulary (summit, …) or refers to the sensitive space ("tip", "peak", "corner"). Then we study the possible evolution of the language used by students of two classrooms. The increase of the geometric vocabulary presented in Figure 4 is obvious, but not similarly in the two classrooms: in classroom 1, the proportion of non geometrical words is still high in session 3, while it is low in session 3 of classroom 2.

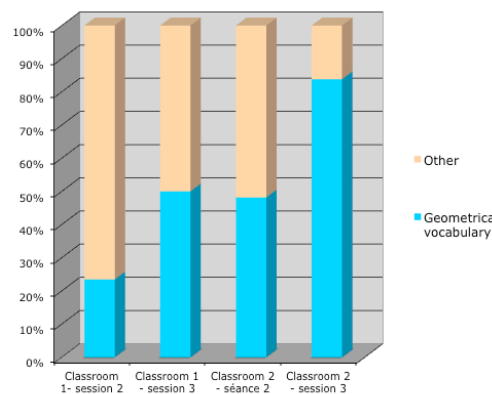

Figure 4: Evolution the proportion of geometrical words between session 2 and 3

In order to explain the differences, we cross these information with the number of geometrical words employed by the teachers (an example is given in Figure 5).

|  | Classroom 1 | Classroom 2 |
|---|---|---|
| Common language used by students | 29 | 3 |
| Geometrical vocabulary used by students | 1 | 6 |
| *Geometrical vocabulary used by teachers* | *2* | *13* |

Figure 5: Vocabulary used by students and teachers during the reproduction of the same model

We find that the role of the teacher is very important in this evolution: teacher 1 doesn't often use the appropriate vocabulary, and even repeats the students' words, while teacher 2 uses systematically a geometrical vocabulary:

Teacher 1: (…) You fold up this tip?//Is this a tip? < *the teacher refers to the center of the PLIOX*>

Students : No!

Teacher 1 : What is it so?

Estelle : A little point.

Teacher 1 : It is not a point. What is it? (…)

Juliette : It's/a trace // when we fold// then/ it leaves marks (…)

The transcriptions highlight the knowledge that appear during the various phases: identification of square, rectangle, triangle and their properties (number of sides and summits, equality of lengths); relationships between square, rectangle, isosceles right-angle triangle (related to the division of the square in two parts); geometrical vocabulary; sides and specific rights of the square (diagonal, median). In addition, according to the transcriptions, we identify the knowledge, which the teachers sometimes grasp (mostly during phase 5): they only concern the identification of square, rectangle, triangle, the geometrical vocabulary (as we have already seen), and the diagonals of the square. Hence we can say that teacher 1 is in acting and not in saying (consequently, he does not propose any institutionalization), while the only objective of teacher 2 is related to the vocabulary (he does not exhibit either relationships between geometrical figures nor their properties).

**CONCLUSION**

The various analysis show the wealth of the PLIOX situation about the way the figures are received, about geometrical and spatial knowledge and the associated vocabulary, sometimes independently of the conditions of implementation. But they also highlight a "dysfunction of a process of teaching" (Margolinas & Laparra, 2008): the mathematical potentialities are not all used in implementations, particularly when the teachers do not identify the issues of learning …